\documentclass[10pt]{amsart}
\usepackage{amsmath, amscd, amsthm, amssymb, amsfonts, verbatim}
\usepackage[dvips]{graphicx,color}
\usepackage[all]{xy}

\title{A dual point of view on the ribbon graph decomposition of moduli space}
\author{Kevin Costello}
\address{Department of Mathematics\\
University of Chicago} \email{costello@math.uchicago.edu}
\date{}

\newcommand{\tr}{\triangle}

\newcommand{\iso}{\cong}
\newcommand{\C}{\mathbb C}

\newcommand{\Z}{\mathbb Z}
\newcommand{\defeq}{\overset{\text{def}}{=}}
\newcommand{\into}{\hookrightarrow}

\newcommand{\op}{\operatorname}

\newcommand{\mbb}{\mathbb}
\newcommand{\mf}{\mathfrak}
\newcommand{\mc}{\mathcal}

\newcommand{\abs}[1]{\left| #1 \right|}

\newcommand{\R}{\mbb R}
\newcommand{\til}{\widetilde}

\newcommand{\opmod}{\overline { \mc N}}
\newcommand{\cmod}{\overline{\mc M}}

\newcommand{\br}{\overline}

\newtheorem{theorem}{Theorem}[subsection]
\newtheorem{proposition}[theorem]{Proposition}

\newtheorem{lemma}[theorem]{Lemma}

\numberwithin{equation}{subsection}

\begin{document}
\maketitle
\section{Introduction}

The ribbon graph decomposition of moduli space is a non-compact
orbi-cell complex homeomorphic to $\mc M_{g,n} \times \R_{> 0}^n$.
This was introduced by Harer-Mumford-Thurston \cite{har1986} and
Penner \cite{pen1987}, and used to great effect by Kontsevich
\cite{kon1992,kon1994} in his proof of the Witten conjecture and his
construction of classes in moduli spaces associated to $A_\infty$
algebras.

In this note,  I discuss in some detail the dual version of the
ribbon graph decomposition of the moduli spaces of Riemann surfaces
with boundary and marked points, which I introduced in
\cite{cos_2004}, and used in \cite{cos_2004oc} to construct
open-closed topological conformal field theories.  This dual version
of the ribbon graph decomposition is a compact orbi-cell complex
with a natural weak homotopy equivalence to the moduli space.  In
the case when all of the marked points are on the boundary of the
surface, the combinatorics of the cell complex is captured by ribbon
graphs, as usual.  In the general case, we find a variant of the
ribbon graph complex.

The idea of the construction is as follows. We use certain partial
compactifications $\opmod_{g,h,r,s}$ of the moduli spaces $\mc
N_{g,h,r,s}$ of Riemann surfaces of genus $g$ with $h>0$ boundary
components, $r$ boundary marked points, and $s$ internal marked
points. The partial compactifications we use are closely related to
the Deligne-Mumford spaces; we allow Riemann surfaces with a certain
kind of singularity, namely nodes on the boundary. The moduli space
$\opmod_{g,h,r,s}$ is an orbifold with corners. The boundary is the
locus of singular surfaces. Therefore the inclusion $\mc N_{g,h,r,s}
\into \opmod_{g,h,r,s}$ is a weak homotopy equivalence of
orbispaces. Inside $\opmod_{g,h,r,s}$ is a natural orbi-cell complex
$D_{g,h,r,s}$, which is the locus where all the irreducible
components of the surface are discs (with at most one internal
marked point). We show that the map $D_{g,h,r,s} \into
\opmod_{g,h,r,s}$ is a weak homotopy equivalence. This shows that
$D_{g,h,r,s} \simeq \mc N_{g,h,r,s}$, giving the desired cellular
model for $\mc N_{g,h,r,s}$.  When $s = 0$, the combinatorics of the
cell complex $D_{g,h,r,0}$ is governed by standard ribbon graphs.
When $s > 0$, we find a variant type of ribbon graph, which has two
types of vertex.

In the standard approach, the non-compact orbi-cell decomposition of
moduli space gives a chain model for the Borel-Moore homology, or
equivalently the cohomology, of moduli space.  The chains are given
by ribbon graphs, and the differential is given by summing over ways
of contracting an edge, to amalgamate two distinct vertices.  The
approach used here gives a compact orbi-cell complex, which gives a
chain model for homology of moduli space.  In the case when $s = 0$,
this is precisely the dual of the standard ribbon graph complex.
That is, the chains are given by ribbon graphs, and the differential
is given by summing over ways of splitting a vertex into two.

The main disadvantage of the approach described here,  compared to
the more traditional approach, is that we do not find a cell complex
homeomorphic to moduli space, but instead a space homotopy
equivalent.

On the other hand, one advantage of this approach is that the cell
complex we find is manifestly compatible with the gluing maps
between the moduli spaces we use. We use ``open-string'' type
gluing; instead of gluing boundary components of surfaces, we glue
together marked points on the boundary.  This allows us to give
\cite{cos_2004,cos_2004oc} a generators-and-relations description
(up to homotopy) for the moduli spaces, considered as a PROP.

In the usual approach to the ribbon graph decomposition, Strebel
differentials are used as a black box to construct the ribbon graph
associated to a surface. Here, ribbon graphs arise in a more
geometric way, as the combinatorial data describing certain very
degenerate surfaces.  The ribbon graph cell complex $D_{g,h,r,s}$,
together with the map (in the homotopy category) $D_{g,h,r,s} \to
\mc N_{g,h,r,s}$, arises immediately from the geometry of the moduli
spaces $\opmod_{g,h,r,s}$. The only part of the construction that
requires any work is showing that the inclusion $D_{g,h,r,s} \into
\opmod_{g,h,r,s}$ is a weak homotopy equivalence. However, the basic
idea of the proof is very simple. The key point is to construct a
deformation retraction of $\opmod_{g,h,0,s}$ onto its boundary. This
is achieved by using the exponential map (for the hyperbolic metric)
to flow the boundary of a surface in $\mc N_{g,h,0,s}$ inwards until
it's singular, which results in a surface in $\partial
\opmod_{g,h,0,s}$.  From considering the fibration $\opmod_{g,h,r,s}
\to \opmod_{g,h,0,s}$ given by forgetting the boundary marked points
and stabilising, we deduce that the inclusion $\partial
\opmod_{g,h,r,s} \into \opmod_{g,h,r,s}$ is also a weak equivalence.
An inductive argument, using the fact that the boundary $\partial
\opmod_{g,h,r,s}$ of $\opmod_{g,h,r,s}$ is a union of products of
similar moduli spaces, allows us to show that the inclusion
$D_{g,h,r,s} \into \partial \opmod_{g,h,r,s}$ is a weak equivalence,
from which the main result follows.

\subsection{Acknowledgements}
I'd like to thank Mohammed Abouzaid for his comments on an earlier
draft of this note.

\section{Moduli of Riemann surfaces with boundary}

\subsection{Riemann surfaces with nodal boundary} A
connected Riemann surface of genus $g$ with $h> 0$  boundary
components has the following equivalent descriptions.
\begin{enumerate}
\item
A compact connected ringed space $\Sigma$, isomorphic as a
topological space to a genus $g$ surface with $h$ boundary
components, and locally isomorphic to $\{ z \in \C \mid \op{Im} z
\ge 0\}$, with its sheaf of holomorphic functions.
\item
A smooth, proper, connected, complex algebraic curve $C$ of genus
$2g-1+h$, with a real structure, such that $C \setminus C(\R)$ has
precisely two components, and $C(\R)$ consists of $h$ disjoint
circles; together with a choice of a component of $C \setminus
C(\R)$.
\item
Suppose $2g - 2 + h > 0$. Then, a Riemann surface with boundary is
equivalently a 2-dimensional connected compact oriented $C^\infty$
manifold $\Sigma$ with boundary, of genus $g$ with $h$ boundary
components, together with a metric of constant curvature $-1$ such
that the boundary is geodesic.
\end{enumerate}
(2) and (3) can be shown to be equivalent (when $2g - 2 + h > 0$) as
follows. Given $\Sigma$,  $C$ is obtained by gluing $\Sigma$ and
$\overline \Sigma$ along their boundary. The real structure on $C$
is arises from the anti-holomorphic involution which is the identity
on $\partial \Sigma$ and interchanges $\Sigma$ and $\overline
\Sigma$.  We denote by $C(\R)$ the set of fixed points of the
anti-holomorphic involution.

Conversely, given $C$, $\Sigma$ is the closure of the chosen
component of $C\setminus C(\R)$ in $C$. The hyperbolic metric on
$\Sigma$ is the restriction of the unique complete hyperbolic metric
on $C$ compatible with the complex structure.

I will also need Riemann surfaces with nodes on the boundary. A
connected Riemann surface with nodal boundary has the following
equivalent descriptions.
\begin{enumerate}
\item
A compact connected ringed space $\Sigma$, locally isomorphic to the
ringed space
$$
\{ (z,w) \in \C \times \C  \mid zw = 0, \op{Im} z \ge 0, \op{Im} w
\ge 0\}
$$
with its sheaf of germs of holomorphic functions. This sheaf is
defined to be the inverse image of the sheaf of germs of holomorphic
functions on $\C \times \C$.
\item
A proper connected complex algebraic curve $C$, with at most nodal
singularities, and a real structure. The real structure on each
connected component $C_0$ of the normalization $\widetilde C$ of $C$
must be of the form $(2)$ above; we also require  a choice of
component of $C_0 \setminus C_0(\R)$. All the nodes of $C$ are
required to be real, that is in $C(\R)$.

We will let $\Sigma \subset C$ be the closure of the chosen
components of $C \setminus C(\R)$.

\item
A compact possibly disconnected Riemann surface with boundary $\til
\Sigma$, together with an unordered finite collection of disjoint
points in $\partial \til \Sigma$, arranged into unordered pairs.
$\Sigma$ is the space obtained from $\til \Sigma$ by identifying
each pair of points on $\partial \til \Sigma$.
\end{enumerate}
As before, to go from the first description to the second, form the
double of the surface $\Sigma$, which is an algebraic curve with a
real structure.

Near a node, $\Sigma$ looks like
$$
\begin{xy}
0*{\Sigma}, 0*\cir<30pt>{r^l}, 0,a(0), **{}, 0+/28pt/*{\bullet},
0+/56pt/*\cir<30pt>{l^r}, 0+/56pt/*{\Sigma}
\end{xy}
$$

The number of boundary components of $\Sigma$ can be defined as
follows. $\partial \Sigma$ will be a  union of circles, glued
together at points as above. Define a smoothing of $\partial
\Sigma$, by replacing each node as above by
$$\xymatrix{
\ar@/_1.2pc/@{-}[rr] &  &  \\
 & \Sigma &  \\
\ar@/^1.2pc/@{-}[rr]&  & }$$ The number of boundary components of
$\Sigma$ is defined to be the number of connected components of this
smoothing.

$\Sigma$ has genus $g$ if it has $h$ boundary components and the
genus of the nodal algebraic curve $C = \Sigma_{\cup \partial
\Sigma} \overline \Sigma$ is $2g-1+h$.

We are also interested in surfaces $\Sigma$ with marked points.
These can be of two types: on the boundary of $\Sigma$, or else in
the interior $\Sigma \setminus \partial \Sigma$. These marked points
must be distinct from the nodes and each other. When $\Sigma$ has
marked points, the double $C$ of $\Sigma$ is an algebraic curve with
marked points, distinct from the nodes. $C$ has a real structure,
and some of the marked points are in $C(\R)$, and some are in
$C\setminus C(\R)$. We say that $\Sigma$ is stable if the double $C$
is, that is if $C$ has only finitely many automorphisms.

Let us suppose that $\Sigma$ is smooth, has non-empty boundary, and
has $(r,s)$ boundary and internal marked points.  Then $\Sigma$ is
unstable if and only if it is a disc and $r + 2s \le 2$ or it is an
annulus and $r = s = 0$. More generally, let $\Sigma$ be a singular
surface. Let $\til \Sigma$ be its normalisation, which is obtained
by pulling apart all the nodes of $\Sigma$.  Each node of $\Sigma$
gives two extra boundary marked points on $\til \Sigma$.  Then
$\Sigma$ is unstable if and only if one of the connected components
of $\til \Sigma$ is.

\subsection{Moduli spaces of surfaces with boundary}

For integers $g,r,s \ge 0,h > 0$,  let $\opmod_{g,h,r,s}$ be the
moduli space of stable Riemann surfaces $\Sigma$ of genus $g$ with
boundary, possibly with nodes on the boundary, with $h$ boundary
components,  with $r$ marked points on the boundary $\partial
\Sigma$, and $s$ marked points in the interior $\Sigma \setminus
\partial \Sigma$. All of the marked points are required to be
distinct from the nodes and each other.

This moduli space is non-empty except for the cases when $g = 0$, $h
= 1$ and $r+2s < 3$, or $g = 0$, $h = 2$ and $r =s = 0$.

Let $\mc N_{g,h,r,s} \subset \opmod_{g,h,r,s}$ be the locus of
non-singular Riemann surfaces (with boundary).

The moduli spaces $\opmod_{g,h,r,s}$ are open subsets of those
constructed by Liu in \cite{liu2002}. Note that in Liu's work,
Riemann surfaces are allowed to have nodes in the interior as well
as on the boundary, whereas the surfaces we use are not allowed to
have nodes in the interior. Also Liu allows the length of boundary
components to shrink to zero, turning boundary components into
punctures. For us, boundary components always have positive length.
Similar moduli spaces were also considered by Fukaya et al.
\cite{fuk_oh_oht_ono2000}. The simplest way to construct these
moduli spaces is to realise that they are very closely related to
the real points of the Deligne-Mumford moduli spaces
$\cmod_{2g-1+h,r+2s}$.
\begin{lemma}
$\opmod_{g,h,r,s}$ is an orbifold with corners of dimension $6g - 6
+ 3h + r+2s$. The interior of $\opmod_{g,h,r,s}$ is $\mc
N_{g,h,r,s}$.
\end{lemma}
Recall that an orbifold with corners is, by definition, an
orbi-space locally modelled on $\R_{\ge 0}^k$.  The reason we get
corners is that there is only one way to smooth any node.   If our
surface has $k$ nodes then it is in a point in the moduli space
locally modelled on the origin in $\R_{\ge 0}^k$.

Let $D_{g,h,r,s} \subset \opmod_{g,h,r,s}$ be the locus consisting
of those surfaces whose irreducible components are all discs with at
most one interior marked point.  By ``irreducible component'' I mean
connected component of the normalisation of the surface.

The main result of this note is
\begin{theorem}
The inclusion $D_{g,h,r,s} \into \opmod_{g,h,r,s}$ is a weak
homotopy equivalence of orbispaces. \label{htpy equiv}
\end{theorem}
The notion of weak homotopy equivalence between orbispaces is
briefly discussed in the appendix.  The proof of this theorem will
be given in the next section.

\begin{figure}
\includegraphics{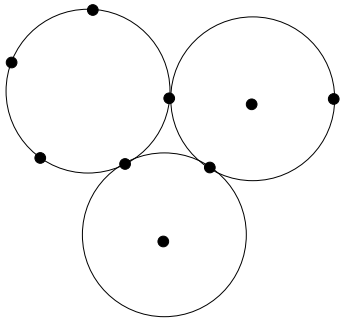}
\caption{A point in $D_{0,2,4,2}$ given by three discs glued
together at the three nodes.}
\end{figure}

\begin{proposition}
Stratify the space $D_{g,h,r,s}$ by saying $\Sigma,\Sigma'$ are in
the same stratum if there exists a homeomorphism $\Sigma \to
\Sigma'$ preserving the marked points and orientation (but not
necessarily respecting the holomorphic structure).

Then this stratification is a decomposition of the compact
orbi-space $D_{g,h,r,s}$ into  orbi-cells. When $r > 0$,
$D_{g,h,r,s}$ is an ordinary space instead of an orbi-space, and the
stratification gives a cell decomposition in the usual sense.
\label{prop cell complex}
\end{proposition}

We will show this by labelling the possible topological types of
surface in $D_{g,h,r,s}$ by a kind of ribbon graph.

Let $\Gamma_{g,h,r,s}$ denote the set of isomorphism classes of
connected  graphs $\gamma$, with the following extra data and
conditions.
\begin{enumerate}
\item
$\gamma$ has $r$ ordered external edges (or tails).
\item
For each vertex $v \in V(\gamma)$, the set of germs of edges at $v$
is cyclically ordered.
\item
The set of vertices $V(\gamma)$ is split into $V_0(\gamma) \amalg
V_1(\gamma)$, and there is given an isomorphism $V_1(\gamma) \iso
\{1,\ldots,s\}$.
\item
All vertices in $V_0(\gamma)$ are at least trivalent, and all
vertices in $V_1(\gamma)$ are at least one valent.
\item
The Euler characteristic of $\gamma$ is $2-2g -h$.
\item
As $\gamma$ is a ribbon graph, we can as usual talk about boundary
components of $\gamma$. There are $h$ unordered boundary components.
\end{enumerate}

Let $\Sigma \in D_{g,h,r,s}$. Associate to $\Sigma$ a graph
$\gamma(\Sigma) \in \Gamma_{g,h,r,s}$.  There is one vertex of
$\gamma(\Sigma)$ for each irreducible component of $\Sigma$, an edge
for each node, and an external edge for each marked point on
$\partial \Sigma$.  As the irreducible components of $\Sigma$ are
all discs, the set of nodes and marked points on the boundary of
each irreducible component has a natural cyclic ordering, coming
from the orientation on the boundary of a disc. Thus the graph
associated to $\Sigma$ has the structure of a ribbon graph. There
are $s$ ordered internal marked points on $\Sigma$, with at most one
on a given irreducible component. A vertex of $\gamma(\Sigma)$ is in
$V_0(\gamma(\Sigma))$ if it doesn't contain an internal marked
point, and it is in $V_1(\gamma(\Sigma))$ if it does.  The
isomorphism $V_1(\gamma(\Sigma)) \iso \{1,\ldots,s\}$ is given by
the ordering of the internal marked points of $\Sigma$.
\begin{lemma}
For any graph $\gamma \in \Gamma_{g,h,r,s}$, the space of $\Sigma
\in D_{g,h,r,s}$ with $\gamma(\Sigma) = \gamma$ is an orbi-cell.
\end{lemma}
\begin{proof}
Let $\gamma \in \Gamma_{g,h,r,s}$, and let $v \in V(\gamma)$. Let
$E(v)$ be the set of germs of edges at $v$.  To give a surface
$\Sigma \in D_{g,h,r,s}$ with $\gamma (\Sigma) = \gamma$ amounts to
giving, for each vertex $v \in V(\gamma)$, a disc $D$ with distinct
points on $\partial D$ labelled by $E(v)$, in a way compatible with
the cyclic order on $E(v)$; and in addition, if $v \in V_1(\gamma)$,
a point in the interior of $D$.  If we change this data by an
automorphism of $\gamma$, then we get an isomorphic surface.

For $v \in V(\gamma)$, let $\til {X}(v)$ be the set of injective
maps $f : E(v) \to S^1$, such that the cyclic order induced on the
image of $f$ by the orientation of $S^1$ coincides with the given
cyclic order on $E(v)$.

If $v \in V_0(\gamma)$ let $X(v) = \til X(v) / PSL_2(\R)$, where
$PSL_2(\R)$ acts on $S^1$ by M\"obius transformations. If $v \in
V_1(\gamma)$ let $X(v) = \til X(v) / S^1$. Note that the spaces
$X(v)$ are cells.

Then the orbi-space of surfaces $\Sigma \in D_{g,h,r,s}$ with
$\gamma(\Sigma) = \gamma$ can be identified with the orbicell
$$
\left( \prod_{v \in V(\gamma)} X(v)\right) / \op{Aut} (\gamma)
$$
where $\op{Aut}(\gamma)$ is the group of automorphisms of $\gamma$
preserving all the labellings.

\end{proof}
This completes the proof of proposition \ref{prop cell complex},
except for the clause about the case when $r > 0$.  This follows
from lemma \ref{lemma no autos} below, which shows that when $r > 0$
surfaces in $\opmod_{g,h,r,s}$ have no non-trivial automorphisms.

One can use this orbi-cell complex to give a chain model for the
rational homology of the moduli spaces $\opmod_{g,h,r,s} \simeq \mc
N_{g,h,r,s}$. A basis for the chain complex is given by ribbon
graphs $\gamma$ in $\Gamma_{g,h,r,s}$ together with an orientation
on the corresponding orbi-cell.  An orientation can be given by
choosing an ordering of the set of vertices of $\gamma$, and at each
vertex an ordering of the set of germs of edges.  It is not
difficult to calculate how changing this ordering changes the
orientation.   The boundary in this chain complex is given by
summing over all ways of splitting a vertex in $V_0$ into two
vertices in $V_0$, and splitting  a vertex in $V_1$ into a vertex in
$V_1$ and a vertex in $V_0$.  In the case $s = 0$, this recovers the
usual ribbon graph model for homology of moduli spaces. When $s >
0$, this chain complex is combinatorially distinct.

When $r > 0$, as $D_{g,h,r,s}$ is in this case an ordinary cell
complex and not an orbi-cell complex, we find a complex computing
the integral homology of moduli space.

Instead of working with an explicit chain complex, I prefer to think
of this result as giving (in the case $s = 0$) a
generators-and-relations description for the PROP controlling open
topological conformal field theory, see \cite{cos_2004,cos_2004oc}.
When $s > 0$ the algebraic statement corresponding to this cell
decomposition is a generators-and-relations description for a
certain natural module over this PROP.

\section{Proof of main theorem}

The idea of the  proof of theorem \ref{htpy equiv} is quite simple,
but a little technical to work out in detail.  The key point is
proposition \ref{prop boundary no points}, which shows that the
orbi-space $\opmod_{g,h,0,s}$ deformation retracts onto its boundary
$\partial \opmod_{g,h,0,s}$,  except for the cases when $(g,h) =
(0,1)$ and $s \in \{0,1\}$.

The first step in the proof is :
\begin{lemma}
There is a map $\opmod_{g,h,r+1,s} \to \opmod_{g,h,r,s}$, which
forgets the last point and stabilises.  This is a locally trivial
fibration in the orbispace sense.
\end{lemma}
\begin{proof}
Recall that there is a map
$$
\cmod_{g,n+1} \to \cmod_{g,n}
$$
$2g-2+n > 0$, given by forgetting the $n+1$'st marked point, and
contracting any resulting unstable components.  This morphism can be
identified with the universal curve on the stack $\cmod_{g,n}$.

There is an induced map
$$
\pi : \opmod_{g,h,r+1,s} \to \opmod_{g,h,r,s}
$$
for $(g,h,r,s)$ stable, that is $4g-4+2h+r +2s > 0$.  This map
removes the $r+1$'st marked point on the surface.  If this leaves
the surface with a disc with two special points, where a special
point is a marked point or a node, then we contract that disc.

We need to show that this is a topologically locally trivial map, in
the orbifold sense.  Let $\Sigma \in \opmod_{g,h,r,s}$, and let us
pick one of the boundary components of $\Sigma$, which we denote by
$\partial_0 \Sigma$. (Recall the definition of boundary components
of a singular surface).  Let us consider adding on a marked point to
$\partial_0 \Sigma$, near a node of $\partial_0 \Sigma$.  If $U
\subset
\partial_0 \Sigma$ is a neighbourhood of a node, then I claim that
the space of ways of adding a marked point in $U$ is just $U$.  It
is clear we can add on a marked point $p$ to a smooth point of $U$.
If $p$ approaches the node,  then we bubble of a disc, which is
inserted into $\partial_0 \Sigma$ at this node. This disc has two
nodes and one marked point, $p$.  There is one way to glue on such a
disc, so we can instead think of adding on a marked point at the
node.   We get the same configuration if $p$ approaches from the
other side.  Therefore the space of possible ways of adding on a
marked point is $U$.

The same behaviour occurs if we looked at a smooth region of a
boundary component. This makes it clear that the map is a locally
trivial fibration.

\end{proof}

\begin{lemma}
If the map $\partial \opmod_{g,h,r,s} \into \opmod_{g,h,r,s}$ is a
weak homotopy equivalence, then so is the map $\partial
\opmod_{g,h,r+1,s} \into \opmod_{g,h,r+1,s}$. \label{prop forget
points}
\end{lemma}
\begin{proof}
Consider the following fibre square.
$$
\xymatrix{
\partial \opmod_{g,h,r+1,s} \ar[d]
\ar[r] & \opmod_{g,h,r+1,s}  \ar[d]
\\
\partial \opmod_{g,h,r,s} \ar[r] & \opmod_{g,h,r,s}}
$$
The vertical arrows are the maps which forget the last marked point.
As the vertical arrows are locally trivial fibrations, and in
particular the fibres are all the same,  the result follows.
\end{proof}

The key step in the proof is the following proposition.

\begin{proposition}
For all $(g,h)$ with $2g-2+h + s > 0$, the inclusion
$$
\partial \opmod_{g,h,0,s} \into \opmod_{g,h,0,s}
$$
is a homotopy equivalence of orbi-spaces. \label{prop boundary no
points}
\end{proposition}

\begin{figure}
\includegraphics{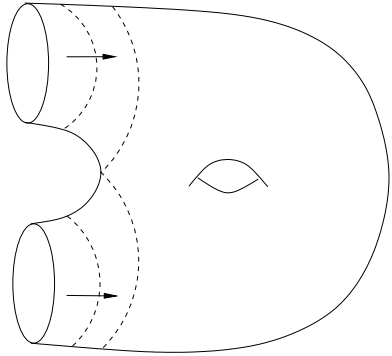}
\caption{By moving the boundaries of the surface inwards, we find a
singular surface.} \label{figure flow}
\end{figure}

\begin{proof}
Let $\Sigma \in \mc N_{g,h,0,s}$.  Note that $\Sigma$ is a smooth
surface. Let $p_1,\ldots,p_s \in \Sigma \setminus \{\partial
\Sigma\}$ be the marked points, and let
$$\Sigma_0 =\Sigma\setminus\{p_1,\ldots,p_s\}.$$ Then
$\Sigma_0$ has a unique hyperbolic metric, compatible with the
conformal structure, such that the boundary is geodesic and such
that the double of $\Sigma_0$ is geodesically complete.  To
construct this metric, observe that the double $C_0$ of $\Sigma_0$
has a unique complete hyperbolic metric, and the anti-holomorphic
involution on $C_0$ is an isometry. It follows that the fixed points
of the involution are geodesic.

Let $V$ be the unit inward pointing normal vector field on $\partial
\Sigma = \partial \Sigma_0$. Using the geodesic flow on $\Sigma_0$,
we can flow $\partial \Sigma$ inwards, as in figure \ref{figure
flow}. That is, for each $s \in \R_{\ge 0}$ in some neighbourhood of
$0$, we have an exponential map
$$
\op{exp}_s : \partial \Sigma  \to \Sigma_0
$$
by flowing the boundary in for a time $s$ along the geodesic flow.
As the double $C_0$ of $\Sigma_0$ is complete, this exists for all
time as a map to $C_0$.  Let $T(\Sigma) \in \R_{> 0}$ be the
smallest number such that the subset
$$
\op{exp}_{T(\Sigma)}( \partial \Sigma  ) \subset \Sigma_0
$$
is singular.
\begin{lemma}
The only singularities of $\op{exp}_{T(\Sigma)} (\partial \Sigma)$
are nodes. In other words, for each $x \in \Sigma_0$,
$\exp_{T(\Sigma)}^{-1}(x)$ consists of at most two points in
$\partial \Sigma$.
\end{lemma}
\begin{proof}
It is easy to convince oneself of the truth of this statement, by
attempting to construct a situation where it fails to hold. It is
nothing to do with hyperbolic geometry, but simply a general fact
about what happens when several moving smooth ``wavefronts'' in two
dimensions collide. At the first collision, at most two can meet at
any point. (It the wavefronts were singular, this would not be
true).

Here is a proof.  Suppose not.  Let $y_1,y_2,y_3 \in \partial
\Sigma$ be three distinct points such that the points
$\exp_{T(\Sigma)}(y_i) \in \Sigma_0$ coincide.  Let $x =
\exp_{T(\Sigma)}(y_i)$.  Let $U_i \subset \partial \Sigma$ be small
neighbourhoods of $y_i$ in $\partial \Sigma$. Note that the
derivative of $\exp_{T(\Sigma)} :
\partial \Sigma \to \Sigma$ never vanishes. We can assume that the
$U_i$ are such that $\exp_{T(\Sigma)} : U_i \to \Sigma$ is
injective.

Consider the three unparameterised $C^\infty$ paths $\gamma_i =
\exp_{T(\Sigma)} (U_i) \subset \Sigma$. These intersect only at $x$.
Suppose $\gamma_1,\gamma_2$ intersect transversely at $x$.  Then,
for some $\epsilon$ sufficiently small, $\exp_{T(\Sigma) -
\epsilon}(U_1)$ and $\exp_{T(\Sigma)-\epsilon}(U_2)$ also intersect
transversely, which contradicts the fact that $T(\Sigma)$ is the
smallest time at which $\exp_{T(\Sigma)}(\partial \Sigma)$ has a
singularity.

Therefore, the tangent vectors to all the paths
$\gamma_1,\gamma_2,\gamma_3$ at $x$ all lie on the same line.  This
implies that their normal vectors do also.  These normal vectors are
non-zero; at least two of them coincide up to rescaling by a
positive real number.  Suppose $\gamma_1,\gamma_2$ have this
property. The geodesic $\exp_{s}(y_i)$ for $0 \le s \le T(\Sigma)$
is normal to $\gamma_i$.  Since $\exp_{s}(y_1)$ and $\exp_{s}(y_2)$
point in the same direction at $x$, it follows that they coincide,
so that $y_1 = y_2$,  which is a contradiction.

\end{proof}

For each $0 \le t \le {T(\Sigma)}$, define a surface
$$
\Sigma(t) \defeq \Sigma \setminus \cup_{s < t} \exp_{s} (\partial
\Sigma)
$$
For $0 \le t < {T(\Sigma)}$, $\Sigma(t)$ is in $\mc N_{g,h,0}$.  The
surface $\Sigma(T(\Sigma))$ is in $\partial \opmod_{g,h,0}$.  To see
this, observe that the previous lemma implies $\Sigma(T(\Sigma))$
has only nodal singularities.  Also there are no unstable
components, simply because there are no hyperbolic polygons with
$\le 2$ sides.

Now define a map of orbispaces
\begin{align*}
\Phi : \mc N_{g,h,0} \times [0,1] & \to  \opmod_{g,h,0} \\
\Phi ( \Sigma, t ) &= \Sigma(t {T(\Sigma)})
\end{align*}
\begin{lemma}
$\Phi$ extends continuously to a map $\opmod_{g,h,0} \times [0,1]
\to \opmod_{g,h,0}$, by defining $\Phi(\Sigma,t) = \Sigma$ for
$\Sigma \in \partial \opmod_{g,h,0}$.
\end{lemma}
\begin{proof}
We have to show that the extension of $\Phi$ so defined is
continuous.   This is not so difficult.  Let $\Sigma_i \in \mc
N_{g,h,0}$ for $i \in \Z_{> 0}$ be a sequence of surfaces converging
to $\Sigma \in \partial \opmod_{g,h,I}$, and let $t_i \in [0,1]$ be
any sequence.    We need to show that
$$\op{lim}_{i \to \infty} \Phi (\Sigma_i, t_i ) = \Sigma$$
Note that  $T(\Sigma_i) \to 0$ as $i \to \infty$. This is well known
in the case of surfaces without boundary. If $C_i \in \mc M_{g,n}$
is a sequence of punctured surfaces converging to a point in
$\partial \br{\mc M}_{g,n}$, then the length of the shortest
geodesic in $C_i$ converges to $0$.   The result in our case follows
from considering the double of the $\Sigma_i$.

Now, since $T(\Sigma_i) \to 0$, it follows that for any sequence
$t_i \in [0,1]$ the sequence
$$\Phi (\Sigma_i, t_i ) = \Sigma_{i}( t_i T(\Sigma_i)) $$
has the same limit as $\Sigma_i$.  That is, $\Phi(\Sigma_i,t_i)$ is
obtained from $\Sigma_i$ by moving the boundary inwards by $t_i
T(\Sigma_i)$, and this tends to zero.  Therefore $\lim
\Phi(\Sigma_i,t_i) = \Sigma$.
\end{proof}
This completes the proof of proposition \ref{prop boundary no
points}. As $\Phi$ gives a deformation retraction of the orbispace
$\opmod_{g,h,0,s}$ onto its boundary $\partial \opmod_{g,h,0,s}$.

\end{proof}

We are nearly finished the proof of theorem \ref{htpy equiv}.
\begin{lemma}
For all $g \ge 0, h \ge 1, r \ge 0, s \ge 0$, with
\begin{align*}
(g,h,r,s) & \neq (0,1,r,0 ) \\
(g,h,r,s) & \neq (0,1,r,1) \\
(g,h,r,s) & \neq (0,2,0,0)
\end{align*}
the inclusion
$$
\partial \opmod_{g,h,r} \into \opmod_{g,h,r}
$$
is a weak homotopy equivalence.
\end{lemma}
\begin{proof}
This follows from \ref{prop boundary no points} and \ref{prop forget
points}, except for the case when $(g,h,r,s) = (0,2,1,0)$. In this
case, it is easy to see that $\opmod_{0,2,1,0} = \R_{\ge 0}$, so
$\partial \opmod_{0,2,1,0} \into \opmod_{0,2,1,0}$ is obviously a
homotopy equivalence.
\end{proof}

\begin{lemma}
Let $\Sigma \in \opmod_{g,h,r,s}$, where $r \ge 1$.  Then $\Sigma$
has no non-trivial automorphisms fixing each of the marked points.
\label{lemma no autos}
\end{lemma}
\begin{proof}
First suppose $\Sigma$ is smooth. Put the hyperbolic metric on
$\Sigma_0 = \Sigma \setminus \{p_1,\ldots,p_s\}$. Any automorphism
of $\Sigma$ preserving the marked points induces an isometry of
$\Sigma_0$, fixing all the marked points on $\partial \Sigma$. Since
the automorphism must act as the identity on the tangent space to
each marked point on $\partial \Sigma$, it must be the identity on a
neighbourhood of each marked point. Since the automorphism is
analytic, it must be the identity everywhere.

Now suppose $\Sigma$ is singular. Let $p \in
\partial \Sigma$ be a marked point. Let $\phi$ be an automorphism of $\Sigma$. Then $\phi$ is
the identity on the irreducible component containing $p$. Suppose
$n$ is a node which joins this component of $\Sigma$ to some other
component. Then $\phi(n) = n$, which implies that $\phi$ is the
identity on the other component at this node. Repeating this
argument we see $\phi$ is the identity everywhere.
\end{proof}

Finally, we can finish the proof of the theorem.
\begin{lemma}
The inclusion $D_{g,h,r,s} \into \opmod_{g,h,r,s}$ is a weak
homotopy equivalence of orbispaces.
\end{lemma}
\begin{proof}
By induction, suppose we have proved the result for all moduli
spaces of lower dimension.

For $k \ge 1$, let $\partial_k \opmod_{g,h,r,s}$ be the space of
surfaces $\Sigma \in \partial \opmod_{g,h,r,s}$, equipped with a map
from the set $\{1,2,\ldots k\}$ to the set of nodes on $\Sigma$.
Lemma \ref{lemma no autos} implies that $\partial_k
\opmod_{g,h,r,s}$ is an ordinary topological space, and not just an
orbi-space.  The spaces $\partial_k \opmod_{g,h,r,s}$ are the $k-1$
simplices of a simplicial space. The face maps are the maps which
forget a node. This simplicial space is the one obtained by iterated
fibre products of the map $\partial_1 \opmod_{g,h,r,s} \to \partial
\opmod_{g,h,r,s}$. Therefore the topological realisation
$\abs{\partial_\ast \opmod_{g,h,r,s}}$ of this simplicial space is
weakly equivalent to $\partial \opmod_{g,h,r,s}$.

Similarly, let $\partial_k D_{g,h,r,s}$ be the space of surfaces in
$D_{g,h,r,s}$ with a map from the set $\{1,2,\ldots k\}$ to the set
of nodes on the surface.  The spaces $\partial_k D_{g,h,r,s}$ form a
simplicial space, and $\abs{\partial_\ast D_{g,h,r,s}}$ is weakly
equivalent to $D_{g,h,r,s}$.

There is a map of simplicial spaces $\partial_\ast D_{g,h,r,s} \to
\partial_\ast \opmod_{g,h,r,s}$.   By induction, we know the maps
$\partial_k D_{g,h,r,s} \to \partial_k \opmod_{g,h,r,s}$ are weak
equivalences.   It follows that the associated map $
\abs{\partial_\ast D_{g,h,r,s}} \to \abs{\partial_\ast
\opmod_{g,h,r,s} } $ on the realisations of our simplicial spaces is
a weak equivalence.

The diagram
$$
\xymatrix{ \abs{\partial_\ast D_{g,h,r,s}} \ar[r]  \ar[d] & \abs{\partial_\ast \opmod_{g,h,r,s} } \ar[d] \\
D_{g,h,r,s} \ar[r] & \partial \opmod_{g,h,r,s} }
$$
commutes, and the vertical arrows and top horizontal arrows are weak
equivalences. It follows that the map $D_{g,h,r,s} \to
\partial \opmod_{g,h,r,s}$ is a weak equivalence, which implies that
$D_{g,h,r,s} \to \opmod_{g,h,r,s}$ is a weak equivalence.

\end{proof}

\section{Appendix : orbispaces} We recall briefly some definitions
from the theory of topological stacks.  See \cite{noo2005} for
details.  We use the word orbispace to refer to a weak topological
Deligne-Mumford stack in the sense of \cite{noo2005}.   An orbispace
is a category fibred in groupoids over the category $\op{Top}$ of
compactly generated Hausdorff topological spaces, satisfying a
descent (or sheaf) condition, and a representability condition. The
Grothendieck topology on the category $\op{Top}$ is that where the
covering maps are the usual open coverings.  The representability
condition is that there exists a surjective map from an ordinary
space which is a local homeomorphism.

Let $X$ be an orbispace, and let $U \to X$ be a surjective local
homeomorphism from a space $U$.  We can form a simplicial space by
taking iterated fibre products of $U$ over $X$. The $n$ simplices
are $(U/X)^{n+1}$, the face maps are projections and the degeneracy
maps are diagonals.  Denote by $\mf{N}^\tr(U/X)$ this simplicial
space, and by $\mf{N}(U/X)$ its geometric realisation. The weak
homotopy type of $\mf{N}(U/X)$ is called the weak homotopy type of
$X$. This is independent of the choice of $U$.

Suppose $f : Y \to X$ is a representable map of orbispaces. This
means that all of the fibres are ordinary spaces. Pick a surjective
local homeomorphism $U \to X$ as above.  Then the map $Y \times_X U
\to Y$ has the same property; in particular $Y \times_X U$ is an
ordinary space.   There is a map $\mf{N}(Y \times_X U / Y )\to
\mf{N}(U/X)$. If this is a weak homotopy equivalence then we say
that the map $f : Y \to X$ is a weak homotopy equivalence.  This
definition can be extended to non-representable maps by refining the
cover $Y \times_X U \to Y$.

To see that this is the correct notion of weak homotopy type of an
orbispace, observe that if $G$ is a discrete group, acting on a
space $X$, and we form the orbispace quotient $X/G$, then the map $X
\to X/G$ is a local homeomorphism, and $\mf{N}(X/ (X/G))$ is one of
the standard models for the homotopy quotient of $X$ by $G$.

\def\cprime{$'$}


\end{document}